\documentclass[12pt]{amsart}

\usepackage{times,epsf,amssymb,amsmath,hyperref, rlepsf, color}

\begin{document}

\newtheorem{thm}{Theorem}[section]
\newtheorem{lem}[thm]{Lemma}
\newtheorem{cor}[thm]{Corollary}

\theoremstyle{definition}
\newtheorem{defn}[thm]{\bf{Definition}}

\theoremstyle{remark}
\newtheorem{rmk}[thm]{Remark}

\def\square{\hfill${\vcenter{\vbox{\hrule height.4pt \hbox{\vrule width.4pt height7pt \kern7pt \vrule width.4pt} \hrule height.4pt}}}$}

\newenvironment{pf}{{\it Proof:}\quad}{\square \vskip 12pt}

\title[Embedded $H$-Planes in $\BH$]{Embedded $H$-Planes in Hyperbolic $3$-Space}
\author{Baris Coskunuzer}
\address{Boston College, Mathematics Department, Chestnut Hill, MA 02467}
\email{coskunuz@bc.edu}
\thanks{The author is partially supported by BAGEP award of the Science Academy.}


\newcommand{\Si}{S^2_{\infty}(\mathbf{H}^3)}
\newcommand{\PI}{\partial_{\infty}}
\newcommand{\BH}{\mathbf{H}^3}
\newcommand{\BR}{\Bbb R}
\newcommand{\BC}{\Bbb C}
\newcommand{\BZ}{\Bbb Z}

\newcommand{\p}{\mathcal{P}}
\newcommand{\e}{\epsilon}
\newcommand{\wh}{\widehat}
\newcommand{\D}{\mathcal{D}}
\newcommand{\s}{\Sigma}
\newcommand{\I}{\mathcal{I}}
\newcommand{\F}{\mathcal{F}}
\newcommand{\A}{\mathcal{A}}
\newcommand{\U}{\mathcal{U}}
\newcommand{\B}{\mathbf{B}}
\newcommand{\C}{\mathcal{C}}
\newcommand{\N}{\mathcal{N}}
\newcommand{\X}{\mathcal{X}}
\newcommand{\cc}{\mathfrak{C}}

\maketitle

\begin{abstract}

We show that for any $\C^0$ Jordan curve $\Gamma$ in $\Si$, there exists an embedded $H$-plane $\p_H$ in $\BH$ with $\PI \p_H =\Gamma$ for any $H\in (-1,1)$. As a corollary, we proved that any quasi-Fuchsian hyperbolic $3$-manifold $M\simeq \Sigma\times\BR$ contains an $H$-surface $\Sigma_H$ in the homotopy class of the core surface $\Sigma$ for any $H\in (-1,1)$. We also proved that for any $C^1$ Jordan curve in $\Si$, there exists a unique minimizing $H$-plane $\p_H$ with $\PI \p_H =\Gamma$ for a generic $H\in(-1,1)$.
\end{abstract}

\section{Introduction}

In this paper, we study the asymptotic $H$-Plateau problem in hyperbolic $3$-space $\BH$. In particular, we finish the existence question for constant mean curvature (CMC) embedded planes in $\BH$ for $H\in(-1,1)$. Throughout the paper, $H$ will represent the mean curvature {\em with sign} where the sign represents the direction of the mean curvature vector with respect to the fixed normal direction in $\BH$ relative to given Jordan curve $\Gamma$ in $\Si$ (See Section \ref{signedMC}). 

In \cite{C1}, we proved that for any Jordan curve $\Gamma$ in $\Si$ with at least one smooth point, there exists an embedded $H$-Plane $\p_H$ in $\BH$ with $\PI \p_H=\Gamma$ for any $H\in (-1,1)$. In that paper, the existence of smooth point condition was crucial to employ the techniques to get the area bounds in compact sets to get convergence. In this paper, by using the results in \cite{C1}, we will remove the smooth point condition from the asymptotic curve $\Gamma$, and finish the asymptotic $H$-Plateau problem in $\BH$ as follows:

\begin{thm} \label{introthm1} Let $\Gamma$ be a $C^0$ simple closed curve in $\Si$. Then, for any $H\in(-1,1)$, there exists a properly embedded, minimizing $H$-plane $\p_H$ in $\BH$ with $\PI \p_H=\Gamma$.
\end{thm}

Note that this result is crucial for the study of quasi-Fuchsian hyperbolic $3$-manifolds. In particular, any quasi-Fuchsian hyperbolic $3$-manifold $M$ naturally induces a limit set $\Lambda_M$ in $\Si$, as $\BH$ is the universal cover of $M$. $\Lambda_M$ is a $\pi_1(M)$-invariant Jordan curve (quasi-circle) with no smooth point in $\Si$ unless it is a round circle (Fuchsian case). Our construction induces canonical $\pi_1$-invariant $H$-planes for a given Jordan curve $\Lambda_M$, and hence, gives the following corollary.

\begin{cor} \label{introcor}
	Let $M$ be a quasi-Fuchsian hyperbolic $3$-manifold, and let $\Sigma$ be the core surface of $M$ with $M\simeq \Sigma\times\BR$. Then, for any $H\in(-1,1)$, there exists an embedded $H$-surface $\Sigma_H$ in $M$ in the homotopy class of $\Sigma$.	
\end{cor}

Notice that a special case for above corollary is the case when $M$ is an almost-Fuchsian hyperbolic $3$-manifold \cite{Uh,HW1}. Almost-Fuchsian hyperbolic $3$-manifolds are a special subfamily of quasi-Fuchsian manifolds, which is "close" to the Fuchsian manifolds. Wang showed that these $H$-surfaces in the corollary give a foliation of $M$ when $M$ is almost-Fuchsian \cite{Wa}.

On the other hand, our techniques give a much stronger result for $C^1$ Jordan curves in $\Si$. In particular, for a $C^1$ Jordan curve $\Gamma$ in $\Si$, the family of $H$-planes $\{\p_H\mid H\in(-1,1)\}$ are pairwise disjoint, and this enables us to show following generic uniqueness result by employing the ideas in \cite{C5}.

\begin{thm} \label{introthm2}
	Let $\Gamma$ be a $C^1$ Jordan curve in $\Si$. Then for a generic $H$ in $(-1,1)$, there exists a unique minimizing $H$-plane $\p_H$ in $\BH$ with $\PI\p_H=\Gamma$.	
\end{thm}

The organization of the paper is as follows: In the next section, we go over the basic notions, and the related results on asymptotic $H$-Plateau problem. In Section 3, we prove Theorem \ref{introthm1} and Corollary \ref{introcor}. In Section 4, we show Theorem \ref{introthm2}. Finally in Section 5, we give some concluding remarks.

\subsection{Acknowledgements} I would like to thank the referee for very valuable remarks and suggestions.


\section{Preliminaries} \label{prelimsec}

In this section, we will go over the basic notions and the related results on asymptotic $H$-Plateau problem. For further details on these notions, and results, one can see the prequel paper \cite{C1} or the survey \cite{C4}.

\subsection{$H$-Plateau Problem} \label{H-Plateausec} \

Constant Mean Curvature (CMC) surfaces are natural generalization of minimal surfaces. While minimal surfaces are critical points of the area functional, CMC surfaces are the critical points of the area functional with a volume constraint.

There are two versions of the area functional with a volume constraint, which are basically equivalent. To describe the parametric version \cite{Hi,Gu}, 
we will take the unit ball $\B$ as the ambient space for simplicity. Note the equivalence of these two versions is discussed in detail in \cite{Gr}. Let $\Gamma$ be a simple closed curve in $\partial\B$. Define the space of immersions $$X_\Gamma=\{u:D^2\to \B \ | \ u\in \C^0(\overline{D^2})\cap \C^2(D^2) \mbox{ and } u(\partial D^2)=\Gamma\}$$

Define the operator $F_H:X_\Gamma\to\BR$ as follows: $$\F_H(u)=\int_{D^2} |u_x|^2+ |u_y|^2 +\dfrac{4}{3}H [u\cdot(u_x\times u_y)] \ dxdy$$

In other words, $\F_H(u)=E(u)+\dfrac{4}{3}H W(u)$ where $$E(u)=\int_{D^2} |u_x|^2+|u_y|^2 \ dxdy \ \mbox{ and } \ W(u)=\int_{D^2} [u\cdot(u_x\times u_y)] dxdy$$
The critical points of $\F_H$ are conformal, $H$-harmonic maps. Hence, for such a $u$, if $\Sigma_u=u(D^2)$ is the image, then $\Sigma_u$ is an immersed $H$-disk. For the Riemannian $3$-manifold setting, see \cite{Gu} and \cite[Section 6]{HK}. In particular for $\BH$ setting, see also \cite[Section 1.2]{Cu}.

There is also a nonparametric formulation of the $H$-Plateau problem in a slightly different way for the extreme curves \cite{C1}. Again, let $\Gamma$ be a Jordan curve in $\partial \B$ and let $\partial \B-\Gamma= E^+\cup E^-$. Let $\Omega^-_u$ be the open region in $\B-\Sigma_u$ which $\Sigma_u\cup E^-$ separates from $\B$, i.e. $\partial \overline{\Omega_u^-}=\Sigma_u\cup E^-$.

For a given $u\in \C^0(\overline{D^2})\cap \C^2(D^2)$, consider the variational problem: $$\I_H(u)=Area(\Sigma_u)+2H.Vol(\Omega^-_u)$$

Here, $Vol(.)$ represents the \underline{oriented} (signed) volume of the region. Notice that for conformal $u$, $E(u)=2 Area(\Sigma_u)$, and $W(u)=3 Vol(\Omega_u^-) +C$ by divergence theorem, where $C$ depends of the choice of the outer surface to define $\Omega_u^-$. Hence, we get $\F_H(u)=2\I_H(u) +C$ which  shows that these two versions are basically equivalent for $H_0$-extreme curves. More precisely, by \cite[Lemma 2.2]{C2},  if $u_0$ minimizes $\F_H$, then $\Sigma_{u_0}$ minimizes $\I_H$. Conversely, if $u_0$ is conformal, and $\Sigma_{u_0}$ minimizes $\I_H$, then $u_0$ minimizes $\F_H$. Note that even though the functionals $\I_H$ and $\F_H$ for $H$-surfaces are equivalent in the first variation, they differ at second variations \cite{Gr}.


\begin{defn} Let $M$ be a $3$-manifold.
\begin{itemize}
\item {\em $H$-surface:} A surface $S$ in $M$ is an $H$-surface if the mean curvature is equal to $H$ everywhere on $S$.


\item {\em Minimizing $H$-disk:} A compact disk $D_u$ in $M$ with $\partial D_u=\Gamma$ is a minimizing $H$-disk in $M$ if $\F_H(u)$ has the smallest value among the disks in $M$ with the boundary $\Gamma$.

\end{itemize}
\end{defn}

Note that when $H=0$, this corresponds to the area minimizing disks. In this case, Meeks and Yau showed that the solution to Plateau problem for  $0$-extreme curves must be embedded \cite{MY1,MY2}. The following lemma is a generalization of their result, showing the existence and the embeddedness of the solutions of $H$-Plateau problem for $H_0$-extreme curves. In particular, we will call a manifold with boundary $H_0$-convex, if along the boundary the mean curvature vector points inside the manifold where $H\geq H_0$ everywhere in $\partial M$. If the boundary is piecewise smooth, we will also require to have the dihedral angles to be less than $\pi$ along the nonsmooth parts of $\partial M$. We will call a curve $H_0$-extreme if the curve is in the boundary of a $H_0$-convex manifold. Here, we will give a simple version of a result proven in \cite{C2}.

\begin{lem} \label{embed} \cite{C2} Let $M$ be a compact $H_0$-convex ball. Let $\Gamma$ be a Jordan curve in $\partial M$. Then, for any $H\in[0,H_0)$, there exists a minimizing $H$-disk $\D_H$ in $M$ with $\partial \D_H=\Gamma$, and any such $\D_H$ is embedded.
\end{lem}

The following lemma is known as maximum principle for $H$-surfaces.

\begin{lem} \cite{Gu} \label{maxlem} [Maximum Principle] Let $\Sigma_1$ and $\Sigma_2$ be two surfaces in a Riemannian $3$-manifold which intersect at a common point tangentially. Let $H_i$ be the (signed) mean curvature of $\Sigma_i$ at the common point with respect to the same normal vector $\mathbf{N}$, i.e. $\mathbf{H}_i=H_i\mathbf{N}$. If $\Sigma_2$ lie in positive side (the normal vector $\mathbf{N}$ direction) of $\Sigma_1$ nearby the common point, then $H_1$ is strictly less than $H_2$, i.e. $H_1 < H_2$.
\end{lem}

\subsection{$H$-Planes in $\BH$:} \label{signedMC} \

Now, we restrict ourselves to $\BH$. Throughout the paper, we will use the notion {\em mean curvature with sign}, and $H\in(-1,1)$. In particular, let $\Gamma$ be a simple closed curve in $\Si$. Fix an orientation on $\Gamma$, and let $\p$ be a plane in $\BH$ with $\PI \p=\Gamma$. Then, the orientation on $\Gamma$ naturally induces an orientation on $\p$, and say the induced normal vector of $\p$ is $\mathbf{N}$. If $\p_H$ is an $H$-plane in $\BH$ with $\PI\p_H=\Gamma$, then the mean curvature vector is $\mathbf{H}=H\mathbf{N}$ where $H\in(-1,1)$ \cite{C3,To}.

\begin{figure}[h]
	\begin{center}
		$\begin{array}{c@{\hspace{.2in}}c}
		
		\relabelbox  {\epsfysize=2in \epsfbox{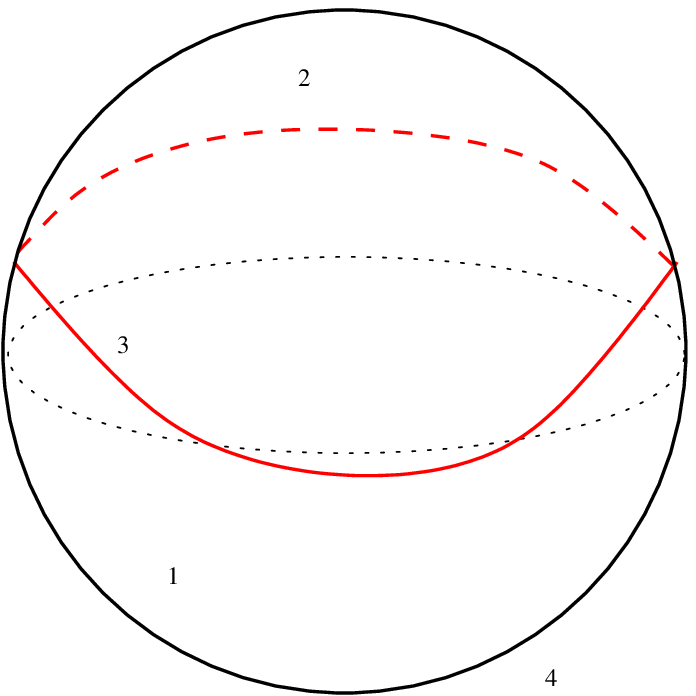}} 
		\relabel{1}{\tiny $\Omega^-$} \relabel{2}{\tiny $\Omega^+$} \relabel{3}{\tiny $\Gamma$} \relabel{4}{\tiny $\Si$} \endrelabelbox &
		
		\relabelbox  {\epsfysize=2in \epsfbox{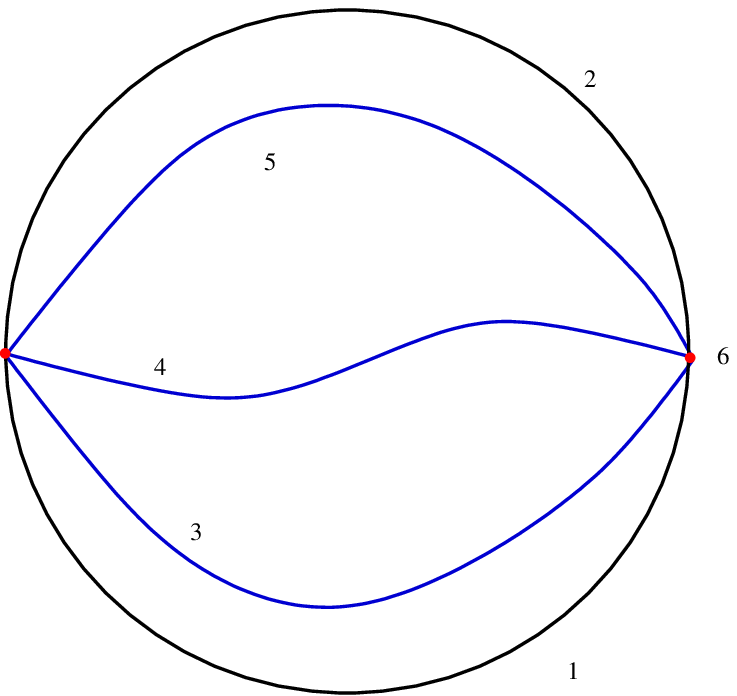}} 
		\relabel{1}{\tiny $\Omega^-$} \relabel{2}{\tiny $\Omega^+$} \relabel{3}{\tiny $\p_{H_2}$} \relabel{4}{\tiny $\p_0$} \relabel{5}{\tiny $\p_{H_1}$} \relabel{6}{\tiny $\Gamma$} \endrelabelbox \\
		\end{array}$
		
	\end{center}
	\caption{\label{fig-signed-H} \footnotesize In the left, $\Si-\Gamma=\Omega^+\cup\Omega^+$. In the right, we depict $\p_{H_1}, \p_0, \p_{H_2}$ (mean curvature with sign) for $H_1<0<H_2$ when $\mathbf{N}$ points upwards.}
\end{figure}

In other words, $\Gamma$ separates $\Si$ into two open disks, say $\Omega^+$ and $\Omega^-$. Fix the normal vector $\mathbf{N}$ on $\p_H$ pointing towards $\Omega^+$ as in Figure \ref{fig-signed-H}. Then the sign of the mean curvature $H$ depends on the direction of the mean curvature vector $\mathbf{H}$, i.e. $H\in(-1,1)$ such that if $H$ is close to $+1$ then $\p_H$ is "close" to $\Omega^-$, and if $H$ is close to $-1$ then  $\p_H$ is "close" to $\Omega^+$. For example, if $\gamma$ is a round circle in $\Si$, then for $K\in(0,1)$ there are two umbilical caps $\p_K^+$ and $\p_K^-$ in $\BH$ with $\PI\p_K^\pm=\gamma$, where $\p_K^\pm$ is an $K$-plane. The mean curvature vectors on $\p_K^+$ and $\p_K^-$ point to each other. With the notation above, we will call $\p_K^+$ as $\p_K$, and we will call $\p_K^-$ as $\p_{-K}$ (See Figure \ref{fig-signed-H}).

Now, by using the definitions above, we define the least area planes and the minimizing $H$-planes in $\BH$.


\begin{defn} [Minimizing $H$-plane] Fix $H\in(-1,1)$. Let $\p_H$ be a complete surface in $\BH$ which is topologically  a disk. We call $\p_H$ a minimizing $H$-plane in $\BH$, if any compact subdisk $D$ in $\p_H$ is a minimizing $H$-disk. Note that the case $H=0$ is special with respect to this definition, and $\p_0$ is known as a {\em least area plane} in the literature.
\end{defn}

For a given complete surface $S$ in $\BH$, we define the {\em asymptotic boundary} of $S$ as follows. If $\overline{\BH}=\BH\cup \Si$ is the natural (geodesic) compactification of $\BH$, and $\overline{S}$ is the closure of $S$ in $\overline{\BH}$, then the asymptotic boundary $\PI S$ of $S$ defined as $\PI S=\overline{S}\cap \Si$.

{\em $H$-shifted convex hulls} $\cc_H(\Gamma)$ are natural generalizations of the convex hulls $\cc(\Gamma)$ in $\BH$ \cite{To, C3}.

\begin{defn} [$H$-shifted Convex Hull] Let $A$ be a subset of $\Si$. Let $\X_A$ be the space of round circles $\gamma$ in $\Si$ with $\gamma\cap A=\emptyset$. Let $D_\gamma$ be the disk which $\gamma$ bounds in $\Si$ with $D_\gamma\supset A$. Fix $H\in (-1,1)$. For each $\gamma\in \X_A$, let $\Delta^H_\gamma$ be the domain in $\BH$ separated by the umbilical $H$-plane (spherical cap) $\p_H$ with $\PI\p_H=\gamma$ and $\PI \Delta^H_\gamma=D_\gamma$. Then, we define the $H$-shifted convex hull as $\cc_H(A)= \bigcap_{\gamma\in \X_A} \Delta^H_A$.
\end{defn}

In particular, when $H=0$, $\cc_0(\Gamma)$ is the usual convex hull of $A$ in $\BH$ which is the smallest convex set in $\BH$ with asymptotic boundary $A$. In other words, $\cc_0(\Gamma)$ is the intersection of all geodesic halfspaces in $\BH$ including $A$ in its asymptotic boundary (supporting halfspaces).

Now, the generalization of the convex hull property of minimal surfaces in $\BH$ to $H$-surfaces in $\BH$ is as follows:

\begin{lem} \label{convexhull} \cite{To, C3}
Let $\Sigma$ be a $H$-surface in $\BH$ where $\PI\Sigma = \Gamma$ and $H\in(-1,1)$. Then $\Sigma$ is in the $H$-shifted convex hull of $\Gamma$, i.e. $\Sigma \subset \cc_H(\Gamma)$.
\end{lem}

We finish this section with the asymptotic regularity result for complete $H$-surfaces in $\BH$ by Tonegawa.

\begin{lem} \label{reglem}\cite{To} [Asymptotic Regularity] 
	Let $\Gamma$ be a $C^1$ Jordan curve in $\Si$. Let $\Sigma_H$ be a complete $H$-surface in $\BH$ with $\PI \Sigma_H=\Gamma$. Then, $\Sigma\cup\Gamma$ is a $C^1$ submanifold with boundary in $\overline{\BH}$.
\end{lem}

Note also that by using some barrier arguments, it is not hard to show that if $\theta_H$ is the intersection angle at infinity between an $H$-surface and the asymptotic boundary $\Si$, then $\cos{\theta_H}=H$ \cite{To}. Note that the compactification of $\BH$ has a natural conformal structure, hence the notion of angle at the asymptotic boundary makes sense.

\section{Existence of Embedded $H$-Planes in $\BH$}

In this section, we prove the main result of the paper. In \cite{C1}, we solved the asymptotic $H$-Plateau problem for Jordan curves with at least one smooth point. In particular, we obtained the following result.

\begin{lem} \label{smoothlem} \cite{C1} Let $\Gamma$ be a simple closed curve in $\Si$ which is differentiable at least at one point. Then, for any $H\in(-1,1)$, there exists a properly embedded, minimizing $H$-plane $\p_H$ in $\BH$ with $\PI \p_H=\Gamma$.
\end{lem}

Even though this result applies to a very general family of curves in $\Si$, the smooth point condition is disturbing, and it prevents a crucial application to low dimensional topology, i.e. quasi-Fuchsian hyperbolic $3$-manifolds. In this section, we will remove this smooth point condition, and finish the asymptotic $H$-Plateau problem in $\BH$.

We will start with the following lemma.

\begin{lem} \label{disjointlem}
	Let $\Gamma_1$ and $\Gamma_2$ be a pair of disjoint $C^1$ Jordan curves in $\Si$. Fix $H\in(-1,1)$. Let $\p_1$ and $\p_2$ be  minimizing $H$-planes in $\BH$ with $\PI\p_i=\Gamma_i$. Then, $\p_1\cap\p_2=\emptyset$.
\end{lem}

\begin{pf} The proof follows from a standard swaping argument. Let $\Si-\Gamma_i=\Omega^+_i\cup\Omega_i^-$. Without loss of generality, assume $\Omega_1^+\subset \Omega_2^+$.
	
Since $\Gamma_i$ is a $C^1$ Jordan curve in $\Si$, $\p_i\cup\Gamma_i$ is $C^1$ submanifold with boundary in $\overline{\BH}$ by Lemma \ref{reglem}. Hence, $\p_i$ separates $\BH$ into two regions, say $\BH-\p_i=\Delta^+_i\cup\Delta^-_i$ where $\PI \Delta^\pm_i = \Omega^\pm_i$. 
	
Assume that $\p_1\cap \p_2\neq \emptyset$. Without loss of generality, assume $H\geq 0$, and the mean curvature vector along $\p_i$ points outside of $\Delta_i^-$.  Hence, by the maximum principle (Lemma \ref{maxlem}), $\p_1\cap \p_2$ cannot be just a collection isolated points, but $\p_1\cap \p_2$ contains a collection of simple closed curves (may not be disjoint), say $\tau$. 

Let $\alpha$ be an innermost curve in $\tau$ where $\alpha$ bounds a disk $D_1$ in $\p_1$ and a disk $D_2$ in $\p_2$. Since $\alpha$ is innermost, $int(D_1)\cap\p_2=int(D_2)\cap\p_1=\emptyset$. Now, recall that $\p_1$ and $\p_2$ are both minimizing $H$-planes in $\BH$, and hence, the disks $D_1$ and $D_2$ are both minimizing $H$-disks with the same boundary $\alpha$, i.e. $I_H(D_1)=I_H(D_2)$.

	
Now, let $T$ be a compact subdisk of $\p_1$ with $D_1\subset int(T)$. Define $T'=(T -D_1)\cup D_2$. Clearly, $T'$ is another embedded disk with $\partial T' = \partial T$. Since $\I_H(D_1)=\I_H(D_2)$, then $\I_H(T)=\I_H(T')$ by construction. As $T$ is minimizing $H$-disk, and $\I_H(T)=\I_H(T')$, then $T'$ is also a minimizing $H$-disk with the same boundary. However, $T'$ has a folding curve along $\alpha$. Therefore, pushing $T'$ to the convex side along $\alpha$  reduces both area and volume, and hence $I_H$. This contradicts to $T'$ being a minimizing $H$-disk. The proof follows.
\end{pf}

Now, we prove the main result of the paper.

\begin{thm} \label{mainthm} [Existence of $H$-Planes] Let $\Gamma$ be a $C^0$ Jordan curve in $\Si$. Let $H\in(-1,1)$. Then, there exists a properly embedded minimizing $H$-plane $\p_H$ in $\BH$ with $\PI \p_H=\Gamma$.
\end{thm}

\begin{pf} Fix $H\in(-1,1)$. $\Gamma$ is a Jordan curve in $\Si$, so it separates $\Si$ into two open disks $\Omega^+$ and $\Omega^-$, i.e. $\Si-\Gamma=\Omega^+\cup\Omega^-$. By using the notation in Section \ref{signedMC}, assume $H\geq 0$ (See Figure \ref{fig-signed-H}).
	
Define a sequence of pairwise disjoint $C^1$ Jordan curves $\{\Gamma_n\}$ in $\Omega^-$ as follows: Let $\Omega_n$ be an open disk in $\Omega^-$ with $\partial \overline{\Omega}_n=\Gamma_n\subset \Omega^-$. We will call $\Gamma_n<\Gamma_m$ if $\overline{\Omega}_n\subset\Omega_m$. We further assume that $\Gamma_n<\Gamma_{n+1}$ for any $n>0$, and $\Omega^-=\bigcup_{n=1}^\infty \Omega_n$. In other words, $\{\Gamma_n\}$ is a {\em monotone} sequence of $C^1$  Jordan curves in $\Omega^-$ converging to $\Gamma$.

By Lemma \ref{smoothlem}, for any $\Gamma_n$, there exists a minimizing $H$-plane $\p_n$ with $\PI\p_n=\Gamma_n$. Note that by Lemma \ref{disjointlem}, $\p_n\cap\p_m=\emptyset$ for any $n\neq m$.

By Lemma \ref{reglem}, $\p_n\cup\Gamma_n$ is a $C^1$ submanifold with boundary in $\overline{\BH}$. Let $\Delta_n$ be the component of 
$\BH-\p_n$ with $\PI\Delta_n=\Omega_n$. Then, by construction, $\overline{\Delta}_n\subset \Delta_{n+1}$ for any $n$. In other words, $\{\Delta_n\}$ is a monotone sequence of open regions in $\BH$ with $\Delta_1\subset\Delta_2\subset ... \subset \Delta_n\subset ...$. Note that  $\Delta_n$ is topologically an open ball in $\BH$ as $\p_n\cup\Gamma_n$ is a closed disk in $\overline{\BH}$.


Let $\Delta^-=\bigcup_{n=1}^\infty \Delta_n$. Then, by construction, $\Delta^-$ is topologically an open ball in $\BH$ with $\PI\Delta^-=\Omega^-$. We claim that $\partial \overline{\Delta^-}=\p_H$ is an embedded, minimizing $H$-plane in $\BH$ with $\PI \p_H=\Gamma$. We first show that $\p_H$ is a smoothly embedded $H$-surface in $\BH$ with $\PI \p_H=\Gamma$. Then, we will finish the proof by showing that $\p_H$ is a plane, i.e. topological type of disk.

First, by convex hull property, Lemma \ref{convexhull}, $\p_H\neq \emptyset$. Now, consider the closure of $\Delta^-$ in the {\em path metric}, say $\wh{\Delta}^-$. In particular, the path metric $\wh{d}(p,q)$ defined as the length of the shortest path in $\Delta^-$ between $p$ and $q$. 

\vspace{.2cm}

\noindent {\bf Claim 1:} Under the path metric, $\wh{\Delta}^-$ is homeomorphic to $\B\cup D$ where $\B$ is the open unit $3$-ball in $\BR^3$, and $D$ is an open disk in $\partial \B$.


\vspace{.2cm}

\noindent {\em Proof of the Claim 1:} Let $\varphi:\B\to \Delta^-$ be a diffeomorphism, and let $g$ be the pull-back metric on $\B$. We consider the closure of $\B$ with the induced topology. First, we show that the region corresponding to $\varphi^{-1}(\Omega^-)$ in the sphere $\partial \B$ will not be in the closure under this metric. Let $p$ be a point in $\Omega^-\subset \Si$, and let $\{p_n\}$ be a sequence in $\Delta^-$ with $p_n \to p$. Let $q_n=\varphi^{-1}(p_n)$. Then by the path metric, the sequence $\{q_n\}$ in $\B$ escapes to infinity, too. This shows that the region corresponding to $\varphi^{-1}(\Omega^-)$ in the sphere $\partial \B$ will not be in the closure under this metric.

On the other hand, we claim that the region corresponding to $\varphi^{-1}(\p_H)$ in the sphere $\partial \B$ will be in the closure under this metric. Similarly, let $p$ be a point in $\p_H$, and let $\{p_n\}$ be a sequence in $\Delta^-$ with $p_n \to p$.
Then, the preimage sequence $\{q_n\}$  will be a Cauchy sequence in $\B$, and so the limit point will be in the closure. In other words, for any point $q\in \partial \B$, and $q_n\to q$, if the sequence $\{p_n\}$ with $p_n=\varphi(q_n)$ converges to a point in compact part of $\BH$ (e.g $p\in\p_H$) then $q\in \partial \B$ will be in the closure. Otherwise, if $\{p_n\}$ escapes to infinity in $\BH$, then so is $\{q_n\}$ in $\B$, and $q\in\partial\B$ will not be in the closure. Hence, the open disk $D$ in $\partial \B$ is exactly the disk $\partial \B-\varphi^{-1}(\Omega^-\cup\Gamma)$.

Now, in order to show that $\wh{\Delta}^-$ homeomorphic to $\B\cup D$, we only need to show that there is no identification in the $D$ under this closure. This is straightforward because of the path metric again. In particular, if $p$ and $q$ are different points in $D$, then the sequences converging to these points are away from each other in $\Delta^-$ with the path metric, even though they might be close in the ambient space $\BH$. With this correspondence, the open ball $\B$ is homeomorphic to $\Delta^-$, and the open disk $D$ in  $\partial\B$ is homeomorphic to $\wh{\p}_H$ where $\wh{\Delta}^-=\Delta^-\cup \wh{\p}_H$. Claim 1 follows.

Hence, $\wh{\p}_H$ is an embedded open disk, which is corresponding to $\p_H$ without possible identifications (self-intersections). In particular, let $\pi:\wh{\p}_H\to\p_H$ be the natural projection map induced by the path metric on $\Delta^-$. Then, $\pi(p)=\pi(q)$ for $p\neq q$ means that there exist sequences $\{p_n\}$ and $\{q_n\}$ in $\Delta^-$ with $p_n\to p$ and $q_n\to q$ where $\wh{d}(p_n,q_m) >\e_0$ for any $m,n$. In other words, $\pi(p)=\pi(q)$ is equivalent to say that in the ambient space $\BH$, $\{p_n\}$ and $\{q_n\}$ converges to the same point, while they are away from each other in $\Delta^-$. 


\vspace{.2cm}

\noindent {\bf Claim 2:} $\wh{\p}_H$ is an $H$-surface.

\vspace{.2cm}

\noindent {\em Proof of the Claim 2:} Let $q$ be a point in $\wh{\p}_H$. Let $B_\e(q)$ be an $\e$-ball around $q$ in $\wh{\Delta}^-$ for sufficiently small $\e>0$ such that $U=\overline{B_\e(q)}$ is topologically a closed ball. Since $\Delta^-$ is exhausted by $\Delta_n$, for sufficiently large $n$, $U\cap \Delta_n\neq \emptyset$. Let $U_n$ be the component of $U\cap \overline{\Delta_n}$ such that $U_n\subset U_{n+1}$ and $U=\bigcup_{n=1}^\infty U_n$. Notice that $E_n=U_n\cap \p_n \subset \p_n\cap B_\e(q)$ is a planar $H$-surface. Since $\p_n$ is a minimizing $H$-plane, and $B_\e(q)$ is a $H'$-convex ball for $H'>>1$, simple area and volume comparison shows that $E_n$ must be a disk. This is because if this is not the case, we can swap the disk in the boundary of $B_\e(q)$ and the disk in $\p_n$ which is not in $B_\e(q)$ with the same boundary, and we can reduce the functional $I_H$. Hence, we obtain a sequence of $H$-disks $\{E_n\}$ with $E_n\subset \p_n$. 

Notice that by construction, $q$ is also a limit point of $\{E_n\}$ as $q\in \partial U$ and $U=\bigcup_{n=1}^\infty U_n$. Since $\{E_n\}$ is a sequence of minimizing $H$-disks, we have a conformal $H$-harmonic parametrization $u_n:D^2\to \BH$ with $u_n(D^2)=E_n$ for any $n$. By composing with a conformal automorphism of the disk, we can also assume $u_n(0)=q_n\to q$. As $u_n$ conformal, the $E(u_n)=2.Area(E_n)$. Then as $\{E_n\}$ are minimizing $H$-disks, there is a uniform bound $K$ such that $E(u_n)<K$ where $E(.)$ is the energy. This is because as we can bound the areas of embedded $H$-disks $\{E_n\}$ in $B_\e(q)$ by the area of $\partial B_\e(q)$, this gives a uniform upper bound $K$ for $E(u_n)$. Then, by using Courant-Lebesque lemma, the equicontinuity of the family $\{u_n\}$ on the compact subsets of $int(D^2)$ follows, see e.g. \cite[Lemma 3.3]{HS}. Hence, by Arzela-Ascoli theorem, there is a continuous function $u:int(D^2)\to \BH$ where $u_n\to u$ pointwise on $int(D^2)$. Furthermore, by using the fact that $\{u_n\}$ are conformal $H$-harmonic maps, again the arguments in \cite{HS} show that $\{u_n\}$ are bounded in the $C^{2,\alpha}$ norm, and $\{u_n\}$ has a subsequence converging $u$ smoothly on $int(D)$. This shows that the disk $E=\wh{\p}_H\cap B_\e(q)$ is a minimizing $H$-disk as the limit of minimizing $H$-disks $\{E_n\}$. This proves that $\wh{\p}_H$ is indeed a smooth $H$-plane. Claim 2 follows. 

\vspace{.2cm}

\noindent {\bf Claim 3:} $\pi:\wh{\p}_H\to\p_H$ is a homeomorphism, and hence $\p_H$ is embedded plane in $\BH$.

\vspace{.2cm}

\noindent {\em Proof of the Claim 3:} Assume that there exist $q_1\neq q_2$ in $\wh{\p}_H$ and $\pi(q_1)=\pi(q_2)=q$. Let $\wh{E}_1$ and $\wh{E}_2$ be small disks in $\wh{\p}_H$ with $q_i\subset \wh{E}_i$. Let $\pi(\wh{E}_i)=E_i\subset \p_H$. Notice that $E_1$ and $E_2$ are both minimizing $H$-disks, and one of them lies in one side of the other one as both lie in $\partial\overline{\Delta^-}$. This contradicts to the maximum principle (Lemma \ref{maxlem}) as $E_1\cap E_2 \supset \{q\}$ and one of them lies in one side of the other one. Indeed, as mean curvature vector points outside $\Delta^-$ along $\p_H$, the mean curvature vectors along $E_1$ and $E_2$ face each other. Note also that for $H>0$, $E_1$ and $E_2$ cannot intersect at an isolated point from concave sides, as the original sequence of embedded planes $\{\p_n\}$ converges to $\p_H$ from the concave side of $\p_H$ ($\{\Delta_n\}$ is a nested sequence). This proves that $\pi:\wh{\p}_H\to\p_H$ is injective. Hence, by construction, $\pi$ is a homeomorphism, and $\p_H$ is an embedded $H$-plane. Claim 3 follows.

\vspace{.2cm}

\noindent {\bf Claim 4:} $\p_H$ is properly embedded in $\BH$.

\vspace{.2cm}

\noindent {\em Proof of the Claim 4:}  Let $\overline{\p}_H$ be the closure  of $\p_H$ in $\BH$. If we can show that $\overline{\p}_H=\p_H$, then this proves that $\p_H$ is properly embedded in $\BH$. Let $q$ be a point in $\overline{\p}_H$. Since $\p_H$ is the limit of $\{\p_n\}$ by construction, then $q$ is in the limit of $\{\p_n\}$, too. In particular, let $\{q^m\}$  be a sequence in $\p_H$ with $q^m\to q$. Then by construction, for each $m>0$, there exists a sequence $\{q^m_n\}$ with $q^m_n\to q^m$ where $q^m_n \in \p_n$. By using diagonal sequence argument, we can conclude that $q$ is a limit point for $\{\p_n\}$ and hence $q$ is indeed in $\p_H$. This proves that $\overline{\p}_H=\p_H$, and hence $\p_H$ is properly embedded in $\BH$. Claim 4 follows.

Now, consider $\overline{\BH}$, the compactification of $\BH$ as a closed ball. Consider the closure of $\Delta^-$ in $\overline{\BH}$, say $\overline{\Delta^-}$. Then, by above, we have $\partial \overline{\Delta^-}=\p_H\cup\Gamma\cup\Omega^-$. This shows that $\PI\p_H=\Gamma$. The proof follows.
\end{pf}

\begin{rmk} [Asymptotic Plateau Problem in $\BH$] Note that for $H=0$, the above theorem shows the existence of a {\em least area plane} for a given $C^0$ Jordan curve in $\Si$ \cite{Ga}. As pointed out in \cite[Remark 3.1]{C1}, there was a subtle problem in the corresponding result in \cite{A2}. With this result, we proved that for a given $C^0$ Jordan curve in $\Si$, there exist a least area plane $\p$ in $\BH$ with $\PI\p=\Gamma$. The main difference between our approach and \cite{A2} is to use $C^1$ curves $\{\Gamma_n\}$ in $\Si$, which bounds properly embedded least area planes $\p_n$. Then, by using the naturally defined regions $\Delta_n$, we converted the problem to codimension-$0$ problem, which bypasses the issues in the approach of \cite{A2}.
\end{rmk}




\subsection{Quasi-Fuchsian Manifolds} \label{QFsec} \

In this section, we will prove the existence of $H$-surfaces in quasi-Fuchsian hyperbolic $3$-manifolds. 

Let $\Sigma$ be a closed, orientable surface of genus $\geq 2$. Let $M=\Sigma\times\BR$ be the product $3$-manifold. Let $g$ be a complete hyperbolic metric on $M$. Then, the universal cover of $(M,g)$ is isometric to $\BH$. Let $\pi:\BH\to M$ be the covering projection. Fix a point $x_0$ in $M$. Let $A=\pi^{-1}(x_0)$ be the set of preimages of $x_0$ in $\BH$. Let $\overline{A}$ be the closure of $A$ in $\overline{\BH}$. Then  $\Lambda_g=\overline{A}\cap\Si$ is called {\em the limit set} of $(M,g)$. If $\Lambda_g$ is a round circle in $\Si$, we call $(M,g)$ a {\em Fuchsian manifold}. If $\Lambda_g$ is a Jordan curve in $\Si$, we call $(M,g)$ a {\em quasi-Fuchsian manifold}. 

The distinction of the family of quasi-Fuchsian manifolds in the space of complete hyperbolic structures on $M$ is coming from $\Lambda_g$ being a Jordan curve. $\Lambda_g$ defines a convex hull $CH(\Gamma)$ in $\BH$, and it is $\pi_1(M)$-invariant. It induces a convex core $\C_M$ in $M$, which simplifies the hyperbolic structure at the ends of $M$. In particular, the ends $M-\C_M$ have a very simple structure as they are foliated by the equidistant surfaces from $\partial \C_M$. By Bers' celebrated theorem, hyperbolic structure on a quasi-Fuchsian manifold $M$ can be uniquely defined by the induced conformal structures at the ends of the manifold, i.e. $\Omega^\pm/\pi_1(M)$ where $\Si-\Gamma=\Omega^+\cup\Omega^-$ \cite{Be}.

Notice that the limit set $\Lambda$ in $\Si$ is $\pi_1$-invariant under the covering transformations. Because of this, the distinction between Fuchsian, and quasi-Fuchsian manifolds are very clear. By group invariance, if $\Lambda_g$ contains one smooth point, then it must be a round circle, and $(M,g)$ would be the Fuchsian manifold. Otherwise, $\Lambda_g$ contains no smooth point, and $(M,g)$ is a quasi-Fuchsian manifold. In this case, $\Lambda_g$ is a nowhere smooth, $C^0$ Jordan curve in $\Si$, and it is called a {\em quasi-circle}.

We should also note that there is another special subfamily of quasi-Fuchsian manifolds called {\em almost-Fuchsian Manifolds}. These are the quasi-Fuchsian manifolds which contain a minimal surface $S$ such that the principle curvatures of $S$ are between $-1$ and $+1$ everywhere. These manifolds are called almost-Fuchsian because their hyperbolic structure is very similar to a Fuchsian manifolds. In particular, the equidistant surfaces $\{S_t\}$ foliate $M$ where $S_t$ is the surface with distance $t\in \BR$ from $S$ \cite{Uh}. Furthermore, $M$ has a unique CMC foliation $\{\mathcal{S}_H\}$ where $H\in(-1,1)$. \cite{Wa}. 

In the following, we will show that the minimizing $H$-planes $\p_H$ constructed in Theorem \ref{mainthm} are canonical, and $\pi_1$-invariant. Hence, we will get an $H$-surface $\Sigma_H$ in $M$ by projecting $\p_H$ in $M$ for $H\in (-1,1)$. In particular, $\Sigma_H$ will be an $H$-surface in the homotopy class of the core surface $\Sigma$.

\begin{cor} \label{QFcor} [Existence of $H$-surfaces in QF-manifolds] 
	Let $M$ be a quasi-Fuchsian hyperbolic $3$-manifold, and let $\Sigma$ be the core surface of $M$ with $M\simeq \Sigma\times\BR$. Then, for any $H\in(-1,1)$, there exists an embedded $H$-surface $\Sigma_H$ in $M$ in the homotopy class of $\Sigma$.	
\end{cor}

\begin{pf} Fix $H\in(-1,1)$. Let $\Lambda$ be the limit set of $M$. Since $\Lambda$ is a $C^0$ Jordan curve in $\Si$, by Theorem \ref{mainthm}, there exists a minimizing $H$-plane $\p_H$ in $\BH$ with $\PI \p_H=\Lambda$. We claim that $\p_H$ is canonical, and $\pi_1$-invariant because of its construction.
	
Let $\Si-\Lambda=\Omega^+\cup\Omega^-$. Recall from the proof of Theorem \ref{mainthm}, $\{\Gamma_n\}$ is a sequence of $C^1$ Jordan curves in $\Omega^-$ with $\Gamma_n\to \Lambda$, and $\p_n$ is a minimizing $H$-plane in $\BH$ with $\PI\p_n=\Gamma_n$ for any $n>0$. Informally, $\p_H$ is the limit of the sequence $\{\p_n\}$.

Let $\beta\in\pi_1(M)$ and $\psi_\beta$ be the corresponding covering transformation. Then, $\psi_\beta$ is an isometry of $\BH$ such that $\psi_\beta(\Lambda)=\Lambda$ and $\psi_\beta(\Omega^\pm)=\Omega^\pm$. We claim that $\psi_\beta(\p_H)=\p_H$. Assume otherwise. Let $\psi_\beta(\p_H)=\p_H'$. As $\psi_\beta$ is an isometry, $\p_H'$ is also a minimizing $H$-plane in $\BH$. Since $\psi_\beta(\Lambda)=\Lambda$, then $\PI \p_H'=\Lambda$, too. 

Let $\Delta^-$ be as in the proof of Theorem \ref{mainthm}, i.e. $\BH-\p_H=\Delta^-\cup\Delta^+$. Let $\psi_\beta(\Gamma_m)=\Gamma_m'$,  $\psi_\beta(\p_m)=\p_m'$, and $\psi_\beta(\Delta_m)=\Delta_m'$.
First assume that $\p_H'\cap\Delta^-\neq \emptyset$. This implies that $\p_{n_o}\cap\p_H'\neq \emptyset$ for sufficiently large $n_o$ as $\Delta^-$ is exhausted by $\{\Delta_n\}$. On the other hand, for sufficiently large $m_0$,  $\Gamma_{n_o}<\Gamma_{m_o}'$ as $\psi_\beta(\Omega^-)=\Omega^-$, and hence $\Gamma_m'\to\Lambda$. Since $\Gamma_{n_o}<\Gamma_{m_o}'$, by Lemma \ref{disjointlem} $\p_{n_o}<\p_{m_o}'$. However, $\p_{m_o}'<\p_H'$ by construction. This implies $\p_{n_o}<\p_{m_o}'<\p_H'$ which contradicts to $\p_{n_o}\cap\p_H'\neq \emptyset$.

Now, we assume $\p_H'\cap\Delta^-= \emptyset$. This implies $\p_H\cap\p_H'=\emptyset$ because of the maximum principle, Lemma \ref{maxlem}. This means $\p_H<\p_H'$, i.e. $\overline{\Delta^-}\subset \psi_\beta(\Delta^-)$. Since $\psi_\beta(\Delta^-)=\bigcup_{m=1}^\infty \Delta_m'$ and $\Gamma_m'<\Gamma$, $\p_H\cap \p_{m_o}'\neq \emptyset$ for sufficiently large $m_0$. Again, as $\psi_\beta(\Omega^-)=\Omega^-$, and $\Gamma_n\to\Lambda$, there exists a sufficiently large $n_0$ such that $\Gamma_{m_o}'<\Gamma_{n_o}$. This means $\p_{m_o}'<\p_{n_o}$ by Lemma \ref{disjointlem}. Hence, we get $\p_{m_o}'<\p_{n_o}<\p_H$ which contradicts to $\p_H\cap \p_m'\neq \emptyset$. 

This shows that for any $\beta\in \pi_1(M)$, $\psi_\beta(\p_H)=\p_H$, and hence $\pi_1$-equivariance of $\p_H$. Then, if $\pi:\BH\to M$ is the covering projection, we get $\pi(\p_H)=\Sigma_H$ an embedded $H$-surface in $M$ in the homotopy class of $\Sigma$. The proof follows.
\end{pf}

\begin{rmk} \label{QFrmk} [Pair of $H$-surfaces]
	Note that the corollary above gives an $H$-surface $\Sigma_H$ for any $H\in(-1,1)$, i.e. mean curvature with sign. If we forget the sign of the mean curvature, for any $K\in(0,1)$, the surfaces $\{\Sigma_K,\Sigma_{-K}\}$ would be a pair of constant mean curvature surfaces in $M$ with mean curvature $K$. 
	
	On the other hand, even in mean curvature with sign case, we might get a pair of $H$-surfaces for some $H\in(-1,1)$ as follows:	
	Notice that in the proof of Corollary \ref{QFcor}, we could have chosen the sequence of $C^1$ Jordan curves $\{\Gamma_n\}$ converging to $\Gamma$ in $\Omega^+$ instead of $\Omega^-$. Then the similar ideas would give us another canonical, $\pi_1$-invariant minimizing $H$-plane $\p_H^+$. If we call the minimizing $H$-plane in the proof above $\p_H^-$, then we end up with two cases. If $\p_H^+=\p_H^-$, then it is not hard to show that $\Gamma$ bounds a unique minimizing $H$-plane in $\BH$ (See \cite{C3} and the proof of Theorem \ref{GUthm} in this paper). If $\p_H^+\neq \p_H^-$, we would have a pair of $H$-surfaces $\Sigma_H^+$ and $\Sigma_H^-$ in $M$ as the covering projection of the planes, i.e. $\pi(\p_H^\pm)=\Sigma_H^\pm$ where $H\in(-1,1)$.	
\end{rmk}

\begin{rmk} \label{AFrmk} [Almost-Fuchsian Manifolds and Uniqueness]
	By using the previous remark, we can conclude another interesting property of almost-Fuchsian manifolds. As mentioned in the introduction of this section, almost-Fuchsian manifolds have an induced CMC foliation \cite{Wa}. Here, each $S_H$ is an $H$-surface for $H\in(-1,1)$. By the maximum principle, this implies $M$ contains a unique $H$-surface for any $H\in(-1,1)$. Hence, by previous remark, if $M$ is almost-Fuchsian, then its limit set $\Lambda$ bounds a unique minimizing $H$-plane in $\BH$ for any $H\in(-1,1)$.	
\end{rmk}

\section{Generic Uniqueness for $C^1$ asymptotic curves} \label{genericsec}

In this section, we will show that if the Jordan curve $\Gamma$ is $C^1$, we have a generic uniqueness result for asymptotic $H$-Plateau problem. In particular, we will show that for a generic $H\in(-1,1)$, $\Gamma$ bounds a unique minimizing $H$-plane in $\BH$. First, we will prove the following lemma by using similar ideas in \cite[Lemma 4.1]{C6}.

\begin{lem} \label{Hdisjointlem}
	Let $\Gamma$ be a $C^1$ Jordan curve in $\Si$. Let $\p_{H_1}$ and $\p_{H_2}$ be minimizing $H_i$-planes with $\PI \p_{H_i}=\Gamma$ where $-1<H_1<H_2<+1$. Then, $\p_{H_1}\cap\p_{H_2}=\emptyset$.	
\end{lem}

\begin{pf} We divide the proof into two steps. First, we show that $\p_{H_1}$ and $\p_{H_2}$ are disjoint near infinity. Then, we will show that they are disjoint in the compact part, too.

\vspace{.2cm} 

\noindent {\bf Step 1:} There exists a sufficiently large ball $\B_{R_0}(p)\subset \BH$ such that $\p_{H_1}-\B_{R_0}(p)$ and $\p_{H_1}-\B_{R_0}(p)$ are disjoint surfaces with the topological type of annulus.

\vspace{.2cm}
	
\noindent {\em Proof of Step 1:} Notice that by Lemma \ref{reglem}, as $\Gamma$ is a $C^1$ Jordan curve in $\Si$, then $\Gamma\cup\p_H$ is a $C^1$ submanifold with boundary in $\overline{\BH}$ where $\p_H$ is the minimizing $H$-plane in $\BH$ with $\PI\p_H=\Gamma$. Furthermore, again by \cite{To}, if $\theta_H$ is the intersection angle at infinity between $\p_H$ and the asymptotic boundary $\Si$, then $\cos{\theta_H}=H$. 
	

Now, let $p\in \BH$ be a point, and $R_0>0$ be sufficiently large that $\partial \B_{R_0}(p)\cap\cc_{H_i}(\Gamma) \subset \{z<\rho\}$ where $\cc_H(\Gamma)$ is the $H$-shifted convex hull. For short, we write $\B_{R_0}$ for $\B_{R_0}(p)$. Let $T_i = \p_{H_i}\cap \B_{R_0}$ and $\gamma_i = \partial T_i$. Then by the previous paragraph, $\gamma_1$ and $\gamma_2$ are disjoint simple closed curves in $\partial \B_{R_0}$ (See Figure \ref{disjointfig}). Again by \cite{To}, we can choose $R_0>0$ sufficiently large so that $\p_{H_i}-\B_{R_0}$ is a graph over $\Gamma\times(0,\rho)$ for some $\rho>0$.  Step 1 follows.

Let $\Delta_i$ be the open region which $T_i$ separates in $\B_{R_0}$ where the mean curvature vector on $T_i$ points outside of $\Delta_i$. Without loss of
generality, we assume $0\leq H_1<H_2<1$. If both negative, just change the direction (See section \ref{signedMC}). If different signs, same argument works as the least area plane $\p_0$ would be a barrier between them. 

Now, we claim that $T_1\cap T_2=\emptyset$. Notice that this will finish the proof as we already know $\p_{H_1} \smallsetminus T_1$ is disjoint from $\p_{H_2} \smallsetminus T_2$ by the first part.

\vspace{.2cm} 

\noindent {\bf Step 2:} $T_1\cap T_2=\emptyset$

\vspace{.2cm} 

\noindent {\em Proof of Step 2:}  Notice that by maximum principle, $T_1\cap T_2$ cannot be just a collection of isolated points. This implies $T_1\cap T_2$ contains a union of closed curves (maynot be disjoint) .

\begin{figure}[b]
	\begin{center}
		$\begin{array}{c@{\hspace{.2in}}c}
		
		\relabelbox  {\epsfysize=1.7in \epsfbox{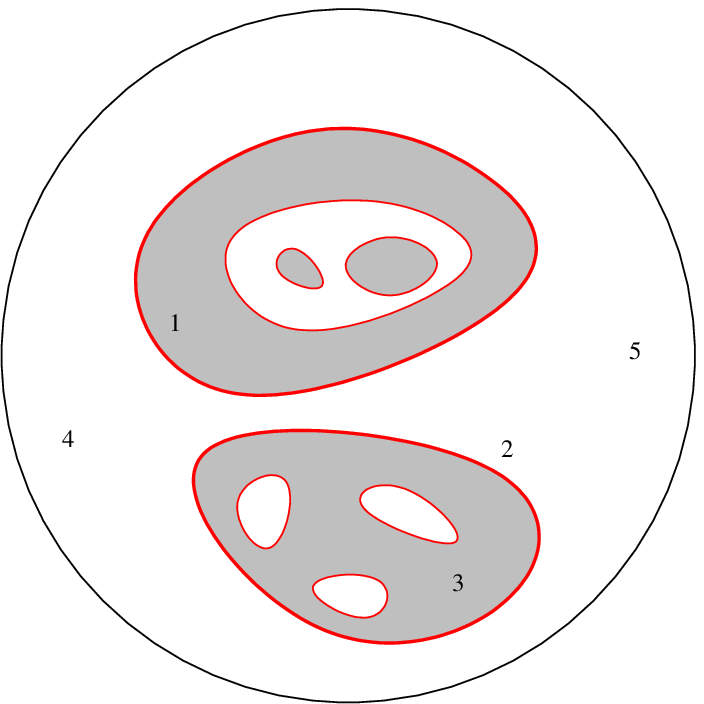}} 
		\relabel{1}{\tiny $\sigma'$} \relabel{2}{\tiny $\alpha$} \relabel{3}{\tiny $\sigma$} \relabel{4}{\tiny $T_1$} \relabel{5}{\tiny $\gamma_1$} \endrelabelbox &
		
		\relabelbox  {\epsfysize=2in \epsfbox{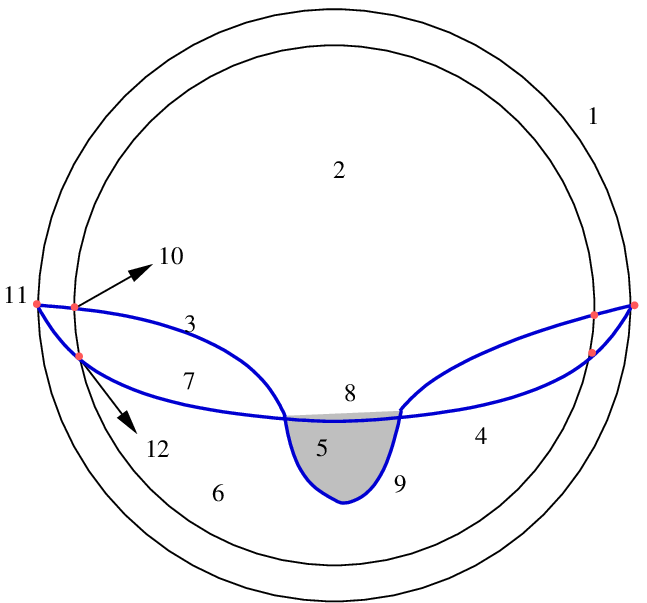}} 
		\relabel{1}{$\Si$}	\relabel{2}{$\B_{R_0}$}	\relabel{3}{$T_1$}
		\relabel{4}{$T_2$}	\relabel{5}{$Q$}	\relabel{6}{$\Delta_2$}
		\relabel{7}{$\Delta_1$}	\relabel{8}{$D_2$}		\relabel{9}{$D_1$} \relabel{10}{$\gamma_1$} \relabel{11}{$\Gamma$} \relabel{12}{$\gamma_2$} \endrelabelbox \\
		\end{array}$
		
	\end{center}
	\caption{\label{disjointfig} \footnotesize In the left, $\sigma$ is an innermost component, and $\alpha$ bounds the disk $D_1 \supset \sigma$.
		In the right, for $0\leq H_1<H_2<1$, $T_1$ is above $T_2$ near the boundary of the ball $\B_{R_0}$. Swaping $D_1$ and $D_2$ gives a contradiction.}
\end{figure}

Consider $T_1\cap\Delta_2$. As $T_1$ is topological type of a disk, $T_1\cap \Delta_2$ is union of planar regions where some of them might be degenerate. In $T_1$, let $\sigma$ be an innermost component of $T_1\cap\Delta_2$. In particular, $\sigma$ does not enclose any other component of $T_1\cap\Delta_2$. For example, in Figure \ref{disjointfig}-left, the shaded regions represent $T_1\cap\Delta_2$. Here, while $\sigma$ is an innermost component in $T_1\cap \Delta_2$, $\sigma'$ is not. Let $\alpha$ be the outer boundary of $\sigma$. Then, $\alpha$ corresponds to a simple closed curve in $T_1\cap T_2$. Then $\alpha$ bounds a disk $D_i$ in $T_i$ where $\sigma\subset D_1$. Let $Q$ be the open region bounded by $D_1\cup D_2$ which may not be connected. Notice that here $D_1\cup D_2$ may not be an embedded sphere, as $int(D_1)\cap int(D_2)$ may not be empty. We fix the orientation on $D_i$ induced by $T_i$, where the mean curvature vector $\mathbf{H}_i$ points outside $\Delta_i$. Then, $D_2-D_1$ is a piecewise smooth sphere with possible self-intersections, whose underlying space is $D_1\cup D_2$. Let the orientation on $Q$ be such that $\partial Q= D_2-D_1$. Here $+$ and $-$ signs represents the corresponding orientations, not the set difference "$\smallsetminus$".





Let $T_2 ' = (T_2-D_2)+D_1$ which is another disk with $\partial T_2'=\gamma_2$ in $\B_{R_0}$. Recall that the signs $\pm$ represents the orientation.  In particular, the underlying set for $T_2'$ is $(T_2\smallsetminus D_2)\cup D_1$ obtained by swaping $D_2$ and $D_1$. Let $T_1 ' = (T_1-D_1) + D_2$ which is another disk in $\B_{R_0}$ with $\partial T_1'=\gamma_1$.  Then, define $\Delta_i'$ be the corresponding region separated by $T_i'$ in $\B_{R_0}$.  

In particular, $\Delta_1'=\Delta_1+Q$, while $\Delta_2'=\Delta_2-Q$. Notice that $\partial \Delta_1'= S_1+T_1'$ and $\partial \Delta_1= S_1 + T_1$ where $S_1$ is the subsurface in the sphere $\partial \B_{R_0}$ separated by $\gamma_1$. Then $\partial (\Delta_1+Q)= \partial \Delta_1+\partial Q=(S_1+T_1)+(D_2-D_1)= S_1+[(T_1-D_1)+D_2]=S_1+T_1'=\partial \Delta_1'$. 

Similarly, $\partial (\Delta_2-Q)= \partial \Delta_2-\partial Q=(S_2+T_2)-(D_2-D_1)=S_2+[(T_2-D_2)+D_1]=S_2+T_2'=\partial \Delta_2'$.
Hence, $\|\Delta_1'\|=\|\Delta_1\|+\|Q\|$ and $\|\Delta_2\|=\|\Delta_2'\|+\|Q\|$ where $\|.\|$ represents oriented (signed) volume (Section \ref{H-Plateausec}).

Since $\p_{H_i}$ is a minimizing ${H_i}$-plane, then $T_i$ is a minimizing $H_i$-disk. Then $\I_{H_i}(T_i)=|T_i| + 2 H_i \|\Delta_i\|$ is minimizer for
the functional $\I_{H_i}$ as defined in the beginning of Section \ref{H-Plateausec}. 

We will get a contradiction by using the fact that $\I_{H_i}(T_i)\leq \I_{H_i}(T_i')$. Since $T_1$ is a minimizing $H_1$-disk, $\I_{H_1}(T_1)\leq \I_{H_1}(T_1')$, i.e. 
$$|T_1| + 2 H_1 \|\Delta_1\| \leq |T_1'| + 2 H_1\|\Delta_1'\|$$ 
This implies 
$$(|D_1|+ |T_1\smallsetminus D_1|) + 2 H_1 \|\Delta_1\| \leq (|D_2|+|T_1\smallsetminus D_1|) + 2 H_1 \|\Delta_1'\|$$
After cancellations, we get 
\begin{equation} \label{H-eqn1}
|D_1| \leq |D_2|+ 2H_1\|Q\|< |D_2|+ 2H_2\|Q\|
\end{equation}

Now, since $T_2$ is a minimizing $H_2$-disk, $\I_{H_2}(T_2)\leq \I_{H_2}(T_2')$, i.e. 
$$|T_2| + 2 H_2 \|\Delta_2\| \leq |T_2'| + 2 H_2\|\Delta_2'\|$$
This implies 
$$(|D_2|+ |T_2\smallsetminus D_2|) + 2 H_2 \|\Delta_2\| \leq (|D_1|+|T_2 \smallsetminus D_2|) + 2 H_2 \|\Delta_2'\|$$. 
After cancellations, we get $|D_2| +2H_2\|Q\| \leq |D_1|$. This contradicts to the strict inequality (\ref{H-eqn1}) above. Hence, Step 2, and the proof of the lemma follows.
\end{pf}




Now, by using the lemma above, we will show that a $C^1$ Jordan curve  in $\Si$ bounds a unique minimizing $H$-plane for a generic $H\in (-1,1)$.

\begin{thm} \label{GUthm}
	Let $\Gamma$ be a $C^1$ Jordan curve in $\Si$. Then for all but countably many $H$ in $(-1,1)$, there exists a unique $H$-plane $\p_H$ in $\BH$ with $\PI\p_H=\Gamma$.	
\end{thm}

\begin{pf} Fix $H\in(-1,1)$. Then, by Theorem \ref{mainthm}, there exists a minimizing $H$-plane $\p_H$ in $\BH$ with $\PI\p_H=\Gamma$. Recall the construction of $\p_H$. Let $\Si-\Gamma=\Omega^+\cup\Omega^-$. We started the construction with a sequence of Jordan curves $\{\Gamma_n^-\}$ in $\Omega^-$ such that $\Gamma_n^-\to \Gamma$. These curves bound minimizing $H$-planes $\p_n^-$, and in the limit, we obtain $\p_n^-\to\p_H^-$. Similarly, start with a sequence of Jordan curves $\{\Gamma_n^+\}$ in $\Omega^+$ with $\Gamma_n^+\to\Gamma$. Let $\p_n^+$ be the minimizing $H$-plane with $\PI \p_n^+=\Gamma_n^+$. Then, similarly we obtain $\p_n^+\to\p_H^+$ where $\p_H^+$ is a minimizing $H$-plane in $\BH$ with $\PI\p_H^+=\Gamma$.
	
\begin{figure}[b]
	\begin{center}
		$\begin{array}{c@{\hspace{.2in}}c}
		
		\relabelbox  {\epsfysize=2in \epsfbox{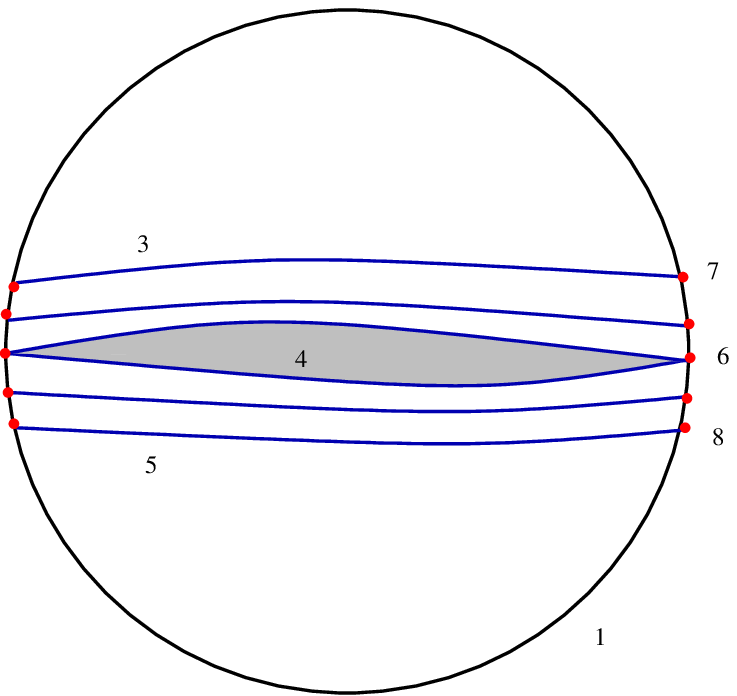}} 
		\relabel{1}{\tiny $\Si$}  \relabel{3}{\tiny $\p_n^+$} 
		\relabel{4}{\tiny $\N_H$} \relabel{5}{\tiny $\p_n^-$} 
		\relabel{6}{\tiny $\Gamma$} \relabel{7}{\tiny $\Gamma_n^+$} \relabel{8}{\tiny $\Gamma_n^-$} \endrelabelbox &
		
		\relabelbox  {\epsfysize=2in \epsfbox{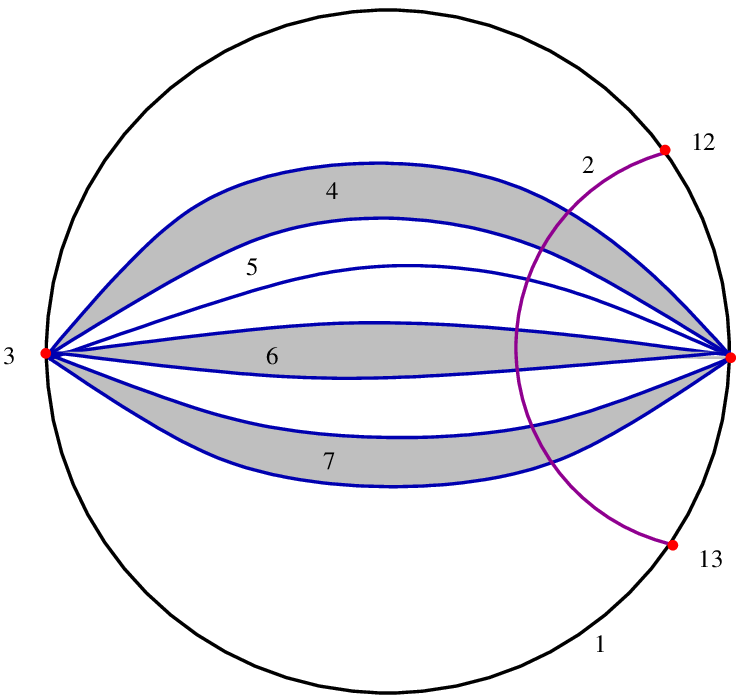}} 
		\relabel{1}{\tiny $\Si$}	\relabel{2}{\tiny $\alpha$}
		\relabel{3}{\tiny $\Gamma$} \relabel{4}{\tiny $\N_{H_4}$}
		\relabel{5}{\tiny $\p_{H_3}$}	\relabel{6}{\tiny $\N_{H_2}$} 
		\relabel{7}{\tiny $\N_{H_1}$}	\relabel{12}{\tiny $q^+$} \relabel{13}{\tiny $q^-$} \endrelabelbox \\
		\end{array}$
		
	\end{center}
	\caption{\label{ngbdfig} \footnotesize In the left, for fixed $H\in(-1,1)$, $\p_n^\pm\to \p_H^\pm$, and defines a canonical region
		$\N_H$ with $\partial \N_H=\p_H^+\cup\p_H^-$.
		In the right, $\N_{H_i}$ is depicted for 4 different values of $H_i\in(-1,1)$ with $H_i<H_{i+1}$. Note that for $H=H_3$, $\N_{H_3}=\p_{H_3}$ is degenerate.}
\end{figure}
	
\vspace{.2cm}
	
\noindent {\bf Claim 1:} If $\p_H^+ =\p_H^-$, then $\Gamma$ bounds a unique minimizing $H$-plane in $\BH$.

\vspace{.2cm}

Assume there exists another minimizing $H$-plane $\Sigma_H$ in $\BH$ with $\PI \Sigma_H=\Gamma$. Let $\p_H=\p_H^\pm$. Let $\BH-\p_H= \Delta^+\cup\Delta^-$.  Then $\Sigma_H\cap \Delta^+\neq \emptyset$ or $\Sigma_H\cap \Delta^-\neq \emptyset$. Say $\Sigma_H\cap \Delta^+\neq \emptyset$. Then, as $\p_n^+\to \p_H$, $\Sigma_H\cap \p_n^+\neq \emptyset$ for sufficiently large $n$. However, since $\Gamma_n^+\cap\Gamma=\emptyset$, by Lemma \ref{disjointlem}, we have $\p_n^+\cap \Sigma_H=\emptyset$ as both are minimizing $H$-planes (See Figure \ref{ngbdfig}-left). This is a contradiction. Claim 1 follows. \ \  $\Box$

\vspace{.2cm}

\noindent {\bf Claim 2:} If $\p_H^+ \neq \p_H^-$, then $\p_H^+\cap\p_H^-=\emptyset$.

\vspace{.2cm}

The proof is very similar to Claim 1. Assume $\p_H^+\cap\p_H^-\neq\emptyset$. Let $\BH-\p_H^+=\Delta^+_+\cup\Delta^-_+$. If $\p_H^+\cap\p_H^-\neq\emptyset$, then $\p_H^-\cap \Delta^+_+\neq \emptyset$. This is because if $\p_H^-\cap \Delta^+_+= \emptyset$ and this would imply $\p_H^-$ and $\p_H^+$ have tangential intersection with one of them lying in one side of the other one where mean curvature vectors point in the same direction. This contradicts to maximum principle, i.e. Lemma \ref{maxlem}. Hence, $\p_H^-\cap \Delta^+_+\neq \emptyset$.

Now, consider $\p_n^+\to \p_H^+$. As $\p_H^-\cap \Delta^+_+\neq \emptyset$, then $\p_H^-\cap \p_n^+\neq \emptyset$ for sufficiently large $n$. However, since $\Gamma_n^+\cap\Gamma=\emptyset$, by Lemma \ref{disjointlem}, we have $\p_n^+\cap \p_H^-=\emptyset$ as both are minimizing $H$-planes. This is a contradiction. Claim 2 follows.\ \  $\Box$


Now, we define a canonical neighborhood $\N_H=[\p_H^-,\p_H^+]$  as follows. If $\p_H^-\neq\p_H^+$, then $\p_H^+\cap\p_H^-=\emptyset$ by Claim 2. Then, let $\N_H$ be the closed region between $\p_H^-$, and $\p_H^+$ such that $\partial \N_H=\p_H^+\cup\p_H^-$ and $\PI \N_H=\Gamma$. 

If $\p_H^+=\p_H^-$, then define $\N_H=\p_H^+$. In this case, we call $\N_H$ a degenerate neighborhood (See Figure \ref{ngbdfig}-right).

\vspace{.2cm}

\noindent {\bf Claim 3:} If $H_1\neq H_2$, then $\N_{H_1}\cap \N_{H_2}=\emptyset$.

\vspace{.2cm}

Assume $H_1<H_2$. Then, by Lemma \ref{Hdisjointlem}, $\p_{H_1}^\pm\cap \p_{H_2}^\pm=\emptyset$ as $\p_{H_i}^\pm$ is a minimizing $H_i$-plane. Then, by definition of $\N_{H_i}$, if $\N_{H_1}\cap \N_{H_2}\neq \emptyset$, then $\p_{H_1}^+\cap \p_{H_2}^-\neq \emptyset$. This is a contradiction. Claim 3 follows. \ \ $\Box$

\vspace{.2cm}

Now, we define the {\em thickness} of the neighborhood $\N_H$. Let $q^\pm$ be a point in $\Omega^\pm\subset \Si$. Let $\alpha$ be a geodesic in $\BH$ with $\PI \alpha =\{q^+, q^-\}$. Notice that $q^+$ and $q^-$ are in the opposite sides of $\Gamma$ in $\Si$ (See Figure \ref{ngbdfig}-right). As $\{\p_H^\pm\}$ are properly embedded, we cans slightly modify $\alpha$ to make it transverse to $\{\p_H^\pm\}$ if necessary. Let $\alpha_H=\alpha\cap \N_H$. Hence, we have $\alpha= \bigcup_{H\in(-1,1)} \alpha_H$. Let $s_H=|\alpha_H|$, the length of $\alpha_H$. We will call $s_H$ the thickness of $\N_H$. Note that this thickness definition depends heavily on the choice of the transverse curve $\alpha$.

Notice that if $s_H=0$, then $\Gamma$ bounds a unique minimizing $H$-plane. Conversely, if $\Gamma$ bounds a unique minimizing $H$-plane, then $s_H=0$ as $\N_H=\p_H^\pm$ and $\alpha_H$ is a point.

For $-1<H_1<H_2<1$, let $\alpha_{[H_1,H_2]}=\bigcup\limits_{H=H_1}^{H_2} \alpha_H$. In particular, $\alpha_{[H_1,H_2]}$ is a finite segment of $\alpha$ between $\N_{H_1}$ and $\N_{H_2}$. In particular, $\alpha_{[H_1,H_2]}$ is the segment of $\alpha$ between the points $p^+,p^-\in \alpha$ where $p^+\in \alpha\cap \p_{H_2}^+$ and $p^-\in \alpha\cap \p_{H_1}^-$. Notice that as $\N_{H}\cap\N_{H'}=\emptyset$ for $H\neq H'$ by Claim 3, then $\alpha_{H}\cap\alpha_{H'}=\emptyset$. Then, the length $|\alpha_{[H_1,H_2]}|=\sum\limits_{H=H_1}^{H_2} s_H$, the sum of $s_H$ for $H\in[H_1,H_2]$. Since the sum $\sum\limits_{H=H_1}^{H_2} s_H$ is finite, this implies $s_H=0$ for all but countably many  $H\in[H_1,H_2]$. This means for all but countably many $H\in[H_1,H_2]$, $\Gamma$ bounds a unique minimizing $H$-plane.

Let $\{t_n \mid n\in\BZ\}$ be a monotone sequence ($t_n<t_{n+1}$) such that $t_n\nearrow +1$ as $n\to\infty$ and $t_n\searrow -1$ as $n\to-\infty$. Hence, $(-1,1)=\bigcup_{-\infty}^\infty [t_n,t_{n+1}]$. If we repeat the same argument, for $\alpha_{[t_n,t_{n+1}]}$, we see that for all but countably many $H\in[t_n,t_{n+1}]$, $\Gamma$ bounds a unique minimizing $H$-plane. Since countable union of countable sets are countable, for all but countably many $H\in(-1,1)$, $\Gamma$ bounds a unique minimizing $H$-plane in $\BH$. The proof follows.	
\end{pf}

\begin{rmk} [Non-minimizing $H$-planes] \label{nonminrem}
	Note that this result is about the generic uniqueness of {\em minimizing $H$-planes}. As minimizing property of $\{\p_H^\pm\}$ is crucial to prove  disjointness of $\N_H$ for different $H$, this result does not imply the uniqueness of "non-minimizing" $H$-planes for a generic $H$ in $(-1,1)$. 	
\end{rmk}

\section{Final Remarks}

\subsection{Asymptotic $H$-Plateau Problem for $C^0$ curves in $\Si$} \

As mentioned in the introduction, in \cite{C1}, we have already proved the existence of embedded $H$-planes for a given Jordan curve $\Gamma$ in $\Si$ when $\Gamma$ has at least one smooth point. Even though this family of Jordan curves in $\Si$ is quite large, there is a very crucial family of Jordan curves, which is left outside: The limit sets of quasi-Fuchsian manifolds. As described in the section \ref{QFsec}, this is an important family of hyperbolic $3$-manifolds, and they induce a natural Jordan curve $\Lambda$ in $\Si$ which is nowhere smooth. Hence, our results in this paper are critical to have applications to quasi-Fuchsian manifolds. 

On the other hand, this result finally finishes off the asymptotic $H$-Plateau problem. When doing this, the most important case $H=0$ has also been resolved. As mentioned in \cite[Remark 3.1]{C1}, there was a subtle issue in the proof of existence of least area planes in $\BH$ for a given Jordan curve $\Gamma$ in $\Si$ \cite{A2}. With our result, this problem has been taken care of, and the asympotic Plateau problem for least area planes was proven, too. Furthermore, for $H=0$ case, Corollary \ref{QFcor} gives an alternative proof for the existence of least area surfaces in quasi-Fuchsian $3$-manifolds.

\subsection{Generalization of Generic Uniqueness Result to $C^0$ curves} \

In Section \ref{genericsec}, we showed that for a given $C^1$ Jordan curve $\Gamma$ in $\Si$, there exists a unique minimizing $H$-plane in $\BH$ for a generic $H\in(-1,1)$. The natural question is that {\em "Can we relax $C^1$ condition on $\Gamma$?} In other words, what if we take a $C^0$ Jordan curve $\Gamma$ in $\Si$? Can we still get uniqueness for minimizing $H$-planes bounding $\Gamma$ for generic $H$? 

The reason we need the $C^1$ condition for our generic uniqueness result is that it is essential for the proof of Lemma \ref{Hdisjointlem}. In particular, for a $C^1$ Jordan curve, we have asymptotic regularity by \cite{To}, which also implies that any $\p_{H_1}$ and $\p_{H_2}$ are disjoint near $\Si$. This is very crucial as $\PI\p_{H_i}=\Gamma$ and they might intersect in an infinite line. Hence, \cite{To} enables us to reduce the problem to the compact case, as it shows the disjointness near infinity. Hence, in the compact part, a straightforward swaping argument gives us Lemma \ref{Hdisjointlem}. However, it is still very interesting question whether this generic uniqueness result can be generalized to $C^0$ Jordan curves or not? Note that the examples of quasi-Fuchsian manifolds constructed by Huang and Wang \cite{HW2} containing arbitrarily many minimal surfaces could be an interesting case to study for this question.

Note also the recent uniqueness results on certain quasi-Fuchsian manifolds by Sanders \cite{Sa} and Seppi \cite{Se}. In particular, let $M\simeq \Sigma\times \BR$ be a quasi-Fuchsian manifold, and $\Lambda$ be its limit set in $\Si$. Sanders showed that if the limit set $\Lambda$ has small Hausdorff dimension ($<1+\delta$), then $\Lambda$ bounds a unique minimal surface in $\BH$. Similarly, Seppi proved that if the quasicircle $\Lambda$ is sufficiently close to the round circle in the Teichmuller sense, $M$ is indeed almost Fuchsian, and hence  the uniqueness of the minimal surface follows. Seppi's result is very interesing as it also implies that almost Fuchsian manifolds contain an open neighborhood of the Fuchsian diagonal in the space of complete hyperbolic structures on $M$. 

 


\subsection{$H$-Laminations of $\BH$ and Quasi-Fuchsian Manifolds} \

Notice that Corollary \ref{QFcor} shows that for any $C^0$ Jordan curve $\Gamma$ in $\Si$, we get two families of embedded $H$-planes $\{\p_H^\pm\}$ with $\PI\p_H^\pm=\Gamma$. In particular, let $\Phi_\Gamma^+=\{\p_H^+\mid H\in(-1,1)\}$ and $\Phi_\Gamma^-=\{\p_H^-\mid H\in(-1,1)\}$. Consider $\Phi_\Gamma^+$. By Corollary \ref{QFcor}, $\Phi_\Gamma^+$ is a pairwise disjoint family of embedded planes in $\BH$ with $\PI\p_H^+=\Gamma$ for any $H\in(1,1)$. Notice that $\Phi_\Gamma^+$ may not be a closed subset of $\BH$. Let $\overline{\Phi_\Gamma^+}$ be its closure in $\BH$. By using the techniques in \cite[Theorem 3.3]{C6}, it can be showed that $\overline{\Phi_\Gamma^+}-\Phi_\Gamma^+$ is a collection of embedded $H$-planes $\wh{\p}_H^+$ for some $H\in(-1,1)$. In particular, $\wh{\p}_{H_o}^+$ is the limit of $H_n$-planes where $H_n\nearrow H_o$. Furthermore, $\overline{\Phi_\Gamma^+}$ is still a pairwise disjoint family of $H$-planes, where for some $H\in(-1,1)$, there are two $H$-planes $\p_H^+$ and $\wh{\p}_H^+$. Hence, $\overline{\Phi_\Gamma^+}$ is a natural lamination of $\BH$ by $H$-planes. Similarly, $\overline{\Phi_\Gamma^-}$ induces another $H$-lamination of $\BH$.

On the other hand, we can apply this idea to the limit sets of quasi-Fuchsian manifolds. Let $M\simeq \Sigma\times \BR$ be a quasi-Fuchsian hyperbolic $3$-manifold, and let $\Lambda$ be its limit set in $\Si$. Then, let $\Phi_\Lambda^+$ be the family of $H$-planes as before, and let $\overline{\Phi_\Lambda^+}$ be its closure. By construction, $\overline{\Phi_\Lambda^+}$ is also $\pi_1$-invariant family of $H$-planes, and if we define $\Sigma_H^+=\Pi(\p_H^+)$ where $\Pi:\BH\to M$ is the covering projection, then we obtain an $H$-lamination $\{\Sigma_H\mid H\in(-1,1)\}$ of $M$, where for some $H\in(-1,1)$, there are two $H$-surfaces $\Sigma_H^+$ and $\wh{\Sigma}_H^+$ as defined before. Again, 
$\overline{\Phi_\Lambda^-}$ induces another $H$-lamination of $M$ by $H$-surfaces.

\subsection{$H$-Foliations of Almost Fuchsian Manifolds} \

Finally, we want to mention that almost Fuchsian manifolds give a very special family of nowhere smooth Jordan curves in $\Si$ in terms of the results in this paper. By \cite{Wa,Uh,HW1}, if $\Lambda$ is the limit set of an almost Fuchsian manifold, then $\Lambda$ induces a foliation $\{\p_H\mid H\in(-1,1)\}$ in $\BH$ where $\p_H$ is a minimizing $H$-plane in $\BH$ with $\PI\p_H=\Lambda$. Furthermore, because of the maximum principle, the foliation automatically implies that $\p_H$ is the unique minimizing $H$-plane with $\PI \p_H=\Lambda$ for any $H\in(-1,1)$. In particular, the limit sets of almost Fuchsian manifolds are very special examples of nowhere smooth Jordan curves in $\Si$  bounding unique $H$-plane $\p_H$ in $\BH$ for any $H\in(-1,1)$.

In particular, if we consider this notion in the terms of the previous section, for almost Fuchsian manifolds, $\Phi_\Lambda^+=\Phi_\Lambda^-$ as $\p_H^+=\p_H^-$ for any $H\in(-1,1)$ by uniqueness. Furthermore, the $H$-lamination $\overline{\Phi_\Lambda^+}$ is indeed an $H$-foliation of $\BH$ by $H$-planes, and the family $\{\Sigma_H\mid H\in(-1,1)\}$ gives a foliation of $M$ by $H$-surfaces.

\end{document}